\documentclass{amsart}

\usepackage{amssymb}
\usepackage[all]{xy}

\newtheorem{thm}{Theorem}

\newtheorem{lem}[thm]{Lemma}
\newtheorem{cor}[thm]{Corollary}

\newtheorem{prop}[thm]{Proposition}

\theoremstyle{definition}
\newtheorem{defn}[thm]{Definition}

\newtheorem{say}[thm]{}
\newtheorem{exmp}[thm]{Example}

\newtheorem*{ack}{Acknowledgments}      

\newtheorem{defn-thm}[thm]{Definition--Theorem}  
\newtheorem{defn-lem}[thm]{Definition--Lemma}  

\theoremstyle{remark}

\renewcommand{\o}[0]{{\mathcal O}} 

\newcommand{\qtq}[1]{\quad\mbox{#1}\quad}
\newcommand{\spec}[0]{\operatorname{Spec}}

\newcommand{\supp}[0]{\operatorname{Supp}}    
    
\newcommand{\codim}[0]{\operatorname{codim}}

\newcommand{\Hom}[0]{\operatorname{Hom}}

\newcommand{\tors}[0]{\operatorname{tors}}

\newcommand{\depth}[0]{\operatorname{depth}}

\def\into{\DOTSB\lhook\joinrel\to}

\def\loccoh#1.#2.#3.#4.{H^{#1}_{#2}(#3,#4)}

\DeclareMathAlphabet{\mathchanc}{OT1}{pzc}%
                                {m}{it}

\usepackage[all]{xy}\xyoption{dvips}

\newcommand{\emb}[0]{\operatorname{emb}}

\newcommand{\assc}[0]{\operatorname{\overline{ass}}}

\begin{document}
\bibliographystyle{amsalpha}


\title{Coherence of  local and global hulls}
\author{J\'anos Koll\'ar}

\maketitle

Let $R$ be a Noetherian, normal, integral domain and $M$ a finite $R$-module.
Its {\it reflexive hull} is the $R$-module
$M^{**}:=\Hom_R\bigl(\Hom_R(M,R),R\bigr)$. The map  $\tau:M\to M^{**}$
is characterized by 3 properties
\begin{itemize} 
\item $M^{**}$ is finite and $S_2$,
\item $\ker \tau=\tors M$, the torsion submodule of $M$ and
\item $\tau_P:(M/\tors M)_P\to M^{**}_P$ is an isomorphism for all height 1 primes $P\subset R$.
\end{itemize}
One can see that $M\mapsto M^{**}$ is a left exact functor from
finite $R$-modules to finite, torsion free, $S_2$ $R$-modules and
it is the ``best'' such functor; see \cite[Tag 0AUY]{stacks-project}.

The concept of hull generalizes this to coherent sheaves over Noetherian schemes.
The notion of torsion subsheaf  is problematic in general, especially
if the support of a sheaf $F$ is not pure dimensional. 
It seems best to assume instead that $\ker\tau$  is  the subsheaf of $F$ generated by those 
local sections whose support is nowhere dense in $\supp F$; we denote
it by $\emb F$. Thus $F/\emb F$ has no embedded associated primes and it is the largest  quotient of $F$ that is $S_1$.

\begin{defn} \label{g.l.hull.defn} Let $X$ be a  Noetherian scheme and $F$ a 
coherent sheaf on $X$. The {\it hull} of $F$ is a 
quasi-coherent sheaf $F^{(**)}$ together with a morphism
$\tau:F\to F^{(**)}$  such that 
\begin{enumerate}
\item $F^{(**)}$ is  $S_2$,
\item $\ker \tau=\emb F$,
\item  $\supp F^{(**)}=\supp F$ and 
\item $\tau_x:(F/\emb F)_x\to F^{(**)}_x$ is an isomorphism for all 
codimension $\leq 1$ points of $\supp F$. 
\end{enumerate}
Let  $Z\subset  X$ be a closed subscheme. 
We also define a local  version of the hull; is is the 
``best'' sheaf that one can associate to $F$ without changing $F|_{X\setminus Z}$.
Thus the {\it local hull} of $F$ centered at $Z$ is a 
quasi-coherent sheaf $F_Z^{(**)}$ together with a morphism
$\tau_Z:F\to F_Z^{(**)}$  such that 
\begin{enumerate}
\item[(1')] $\depth_ZF_Z^{(**)}\geq 2$,
\item[(2')] $\ker \tau_Z=\tors_Z F$, the largest subsheaf of $F$ supported on $Z$,
\item[(3')]  $\supp F^{(**)}$ is the closure of $\supp(F|_{X\setminus Z})$ and 
\item[(4')] $\tau_Z:F\to F_Z^{(**)}$ is an isomorphism over $X\setminus Z $.
\end{enumerate}
(The depth for non-coherent sheaves is defined in   \cite[Exp.III]{sga2}.
By that definition (\ref{g.l.hull.defn}.1') is equivalent to  (\ref{g.l.hull.defn}.5'); see Lemma \ref{rel.hull.depth.2.lem}.)

If $x\in \supp F$ is a point, we let $F_x^{(**)}$ denote the hull of the localization
$F_x$ centered at $\{x\}\subset X_x$, and call it  the
 {\it punctual hull} of $F$ at $x$.

It is easy to see that local and global hulls always exists.
In the local case
$$
F_Z^{(**)}=j_*\bigl(F|_{X\setminus Z}\bigr),
\eqno{(\ref{g.l.hull.defn}.5')}
$$
where $j:X\setminus Z\into X$ is the natural injection.
(See Lemma \ref{rel.hull.depth.2.lem} for condition (1').)
In the global case
$$
F^{(**)}=\varinjlim\ (F/\emb F)_Z^{(**)},
\eqno{(\ref{g.l.hull.defn}.5)}
$$
where $Z$ ranges through all codimension $\geq 2$ subsets of
$\supp F$. 
\end{defn}

The formulas (\ref{g.l.hull.defn}.5--5')  show that  local and global hulls are quasi-coherent.
Our aim is to understand when  local and global  hulls are coherent, using
\begin{itemize} 
\item  $\assc(F)$, the set of closures $W_i\subset X$ of the associated points of $F$ and
\item  properties of the localizations $F_x$  for 
all points $x\in Z$.
\end{itemize}

The local case turns out to be the easier.
The following theorem is a  generalization of   \cite[IV.5.11.1]{EGA}.

\begin{thm}\label{pushfd.finite.thm}
Let $X$ be a  Noetherian  scheme,  $Z\subset X$ a  closed subscheme and  $F$  a coherent sheaf on $X$ such that $\tors_ZF=0$. 
The following are equivalent.
\begin{enumerate}
\item The local hull $F^{(**)}_Z$ is coherent. 
\item For every   $x\in Z$, the punctual hull
$F^{(**)}_x$  is coherent.
\item For every   $x\in Z$, $H^1_x(X_x, F_x)$ has finite length  over $k(x)$.
\item For every  $x\in Z$, the completion
$\hat{F}_x$ has no 1-dimensional associated primes.
\item For every $W\in \assc(F)$ the local hull $(\o_W)^{(**)}_{Z\cap W}$ is coherent. 
\item For every $W\in \assc(F)$  and  $x\in Z\cap W$, the punctual hull
$(\o_W)^{(**)}_x$  is coherent.
\item For every $W\in \assc(F)$  and  $x\in Z\cap W$, $H^1_x(W_x, \o_{W_x})$ has finite length  over $k(x)$.
\item For every $W\in \assc(F)$  and  $x\in Z\cap W$, the completion
$\hat{\o}_{x,W}$ has no 1-dimensional associated primes.
\end{enumerate}
\end{thm}

For the  coherence of global hulls, there are some
obvious restrictions. 

\begin{say}[Necessary conditions] \label{purity.cond.lem} 
Let $X$ be a  Noetherian  scheme and  $F$  a coherent, $S_2$ sheaf on $X$.
We prove in (\ref{DEPTH.OF.SUBMODS}) that its support
satisfies the following  {\it codimension 1 purity condition.}
\begin{enumerate}
\item If $x\in \supp F$ has codimension 1 in some  irreducible component  of $\supp F$ then $x$ has codimension 1 in $\supp F$. 
\end{enumerate}
Note that the condition is vacuous if $\supp F$ is irreducible, but it
is meaningful if $\supp F$ is not pure dimensional. For instance,
the union of a plane and an intersecting line can not be the support
of a   coherent, $S_2$ sheaf.
(It is, however, the  support
of a   quasi-coherent, $S_2$ sheaf.)

Since a coherent sheaf $F$ and its hull $F^{(**)}$ have the same support by definition,
if the hull is coherent then $\supp F$ also satisfies condition (1).

Another condition is the following.
If $F^{(**)}$ is coherent then $\tau:(F/\emb F)\to F^{(**)}$ is an
isomorphism over a dense open subset of $\supp F$, hence
\begin{enumerate}\setcounter{enumi}{1}
\item there is a dense open subset $U\subset \supp F$ such that 
$(F/\emb F)|_U$ is $S_2$. 
\end{enumerate}
This condition turns out to be automatic if $X$ is ``nice''
(for instance   excellent or  N-2) but not always; see Example \ref{bad.exmps.2}.
\end{say}

Thus the key question is to understand when condition (\ref{purity.cond.lem}.1)
is also sufficient for the existence of coherent hulls. The answer is given by the following.

\begin{thm} \label{glob.hull.conds.thm}
For  a Noetherian scheme $X$ the following are equivalent.
\begin{enumerate}
\item A coherent sheaf $F$ has a coherent hull iff the codimension 1 purity condition (\ref{purity.cond.lem}.1) holds for $\supp F$.
\item $\o_W$ has a coherent hull for  every integral subscheme $W\subset X$.
\item For every integral subscheme $W\subset X$
\begin{enumerate}
\item there is an open, dense subset $W^0\subset W$ that is $S_2$ and
\item for every point  $x\in W$ of codimension $\geq 2$, the completion
$\hat{\o}_{x,W}$ has no 1-dimensional associated primes.
\end{enumerate}
\end{enumerate}
\end{thm}

\begin{say}\label{loc.glob.tf.exmp}
The relationship  between the global and local hulls is  clear 
in most cases. Assume that $F$ is coherent, $S_1$ and let
 $Z\subset \supp F$ be a closed subscheme  of codimension $\geq 2$. 
Then 
\begin{enumerate}
\item  $F^{(**)}_Z\subset F^{(**)}$,
\item  $F^{(**)}_Z= F^{(**)}$ iff $F|_{X\setminus Z}$ is $S_2$ and 
\item  in general  $F^{(**)}_Z\subset F^{(**)}$ is determined by the exact sequence
$$
0\to F\to F^{(**)}_Z\to \tors_Z\bigl(F^{(**)}/F\bigr)\to 0.
$$
\end{enumerate}
  So there is no conflict between the
local and global notions. However, if 
$\codim_XZ=1$ then the local hull $F^{(**)}_Z$ is neither coherent nor a subsheaf of the hull $F^{(**)}$. 
\end{say}

\begin{defn}[Pure hulls]
If $X$ is a finite type scheme, more generally if $X$ admits a
dimension function (cf.\ \cite[Tag 02I8]{stacks-project}), then one can define another version of the hull---denoted  by
$F^{[**]}$---where the kernel of $F\to F^{[**]}$ is not $\emb F$
but the subsheaf of sections whose support has dimension
$<\dim F$. Thus $F^{[**]}$ has pure dimensional support and so
it could be called the {\it pure hull} of $F$.

This was the definition adopted in \cite{k-hh}
and  \cite{k-modbook} and it is better suited to the applications there.

 The two notions agree if $\supp F$ is irreducible, more generally, if
$\supp F$ is pure dimensional. Thus Theorem \ref{glob.hull.conds.thm}
takes the following form for pure hulls.
\end{defn}

\begin{thm} \label{pure.glob.hull.conds.thm}
Let $X$ be  a Noetherian scheme with a dimension function. The following are equivalent.
\begin{enumerate}
\item Every coherent sheaf $F$ has a coherent pure hull $F^{[**]}$.
\item $\o_W$ has a coherent hull for  every integral subscheme $W\subset X$.
\item For every integral subscheme $W\subset X$
\begin{enumerate}
\item there is an open, dense subset $W^0\subset W$ that is $S_2$ and
\item for every point  $x\in W$ of codimension $\geq 2$, the completion
$\hat{\o}_{x,W}$ has no 1-dimensional associated primes. \qed
\end{enumerate}
\end{enumerate}
\end{thm}

\begin{say}[What is new in this note?]
Most of the results in this note are contained in---or can be obtained by a careful contemplation of---\cite[IV.5.11.1]{EGA}, which essentially says that local hulls are coherent
if the completed local rings $\hat\o_{x,W}$  are reduced and pure dimensional. 

Our observation is that only 1-dimensional associated primes 
of the completed local rings $\hat\o_{x,W}$ cause problems, and this way
one obtains necessary and sufficient conditions.
The key technical point, following \cite{k-normal},  is the
systematic use of punctual hulls.
\end{say}

\begin{ack} I thank A.~J.~de~Jong, W.~Heinzer, R.~Heitmann, S.~Kov\'acs
and S.~Loepp   for helpful
comments and references.
Partial financial support    was provided  by  the NSF under grant number
 DMS-1362960.
\end{ack}

\section{Proof of Theorem \ref{pushfd.finite.thm}}

\begin{say}[Left exactness] \label{left.exact.say}
Since push-forward is a left exact functor, the  formula (\ref{g.l.hull.defn}.5')
shows that the local hull $F\mapsto F^{(**)}_Z$ is also left exact.

The situation is more delicate for the hull. The problem is that,
even if $F_1\into F_2$ is an injection, it can happen that a subscheme
$Z$ has codimension $\geq 2$ in $\supp F_2$ but
codimension $1$ in $\supp F_1$. This is, however, the only obstruction.
The following special case is especially useful.

\medskip

{\it Claim \ref{left.exact.say}.1.} Let $F$ be a coherent,  $S_1$ sheaf
and $G\subset F$ a subsheaf.
Assume that $\supp F$ satisfies the purity condition (\ref{purity.cond.lem}.1). Then
$G^{(**)}\subset F^{(**)}$.\qed
\end{say}

We start the proof of Theorem \ref{pushfd.finite.thm}
with the implications
(\ref{pushfd.finite.thm}.1) $\Leftrightarrow$ (\ref{pushfd.finite.thm}.5).

\begin{say}[Comparing  hulls] \label{tf.pf.say}
 Let $X$ be a  Noetherian, affine, integral  scheme
and $F$  a torsion-free coherent sheaf on $X$. If $F$ has generic rank $r$ then there are injections
$$
\o_X^r\into F\qtq{and} F\into \o_X^r.
\eqno{(\ref{tf.pf.say}.1)}
$$
Thus we get injections
$$
\bigl(\o_X^{(**)}\bigr)^r\into F^{(**)}\qtq{and}  F^{(**)}\into 
\bigl(\o_X^{(**)}\bigr)^r.
\eqno{(\ref{tf.pf.say}.2)}
$$
If $Z\subset X$ is a closed, 
nowhere dense subscheme, we also get injections
$$
\bigl((\o_X)^{(**)}_Z\bigr)^r\into F^{(**)}_Z\qtq{and}  F^{(**)}_Z\into 
\bigl((\o_X)^{(**)}_Z\bigr)^r.
\eqno{(\ref{tf.pf.say}.2')}
$$
Thus we conclude that
\medskip

{\it Claim \ref{tf.pf.say}.3.} $F^{(**)}$ (resp.\ $F^{(**)}_Z$) is coherent iff
$\o_X^{(**)}$ (resp.\ $(\o_X)^{(**)}_Z$)  is. \qed
\medskip

Let us now drop the assumptions that $X$ is integral and $F$ is torsion free. 
Let $W_i\subset X$ be the closures of the associated points of $F$.
\medskip

{\it Claim \ref{tf.pf.say}.4.} Assume that $\supp F$ satisfies the purity condition (\ref{purity.cond.lem}.1). Then
$F^{(**)}$ (resp.\ $F^{(**)}_Z$) is coherent iff each
$\o_{W_i}^{(**)}$ (resp.\ $(\o_{W_i})^{(**)}_{Z\cap W_i}$)  is.
\medskip

Proof. By d\'evissage  (see, for instance, \cite[Tag 01YC]{stacks-project})
$F$ admits a finite filtration 
$0=G_0\subset G_1\cdots\subset G_m=F$ such that every graded piece
$G_{r+1}/G_r$ is torsion free over some $W_i$. By induction we obtain that
if each
$\o_{W_i}^{(**)}$ is coherent then so is   $F^{(**)}$. 
Conversely,  for each $i$ there is an injection
$ \o_{W_i}\into F$. Thus if $F^{(**)}$ is coherent then so is
$\o_{W_i}^{(**)}$. The proof of the local version is the same. \qed
\end{say}

Next we show some cases when the structure sheaf has a coherent hull,
using some  better known conditions of the theory of commutative rings.
Note that (\ref{conds.for.coh.hull}.1) is known as the N-1 condition   
\cite[Tag 0BI1]{stacks-project}. 
If, in addition,  $W$ is universally catenary,
 then  $\pi:W^{\rm n}\to W$ preserves
codimension \cite[Tag 02II]{stacks-project}. In particular,  if $X$ is
 excellent then both of these conditions hold for every integral subscheme $W\subset X$, but the conditions
(\ref{conds.for.coh.hull}.1--2) are strictly weaker than excellence.

\begin{prop} \label{conds.for.coh.hull}
Let $W$ be  a Noetherian, integral scheme such that
\begin{enumerate}
\item the normalization  $\pi:W^{\rm n}\to W$ is finite and
\item if $w'\in W^{\rm n}$ has codimension 1 then $\pi(w')\in W$ also has
codimension 1.
\end{enumerate}
Then the hull $\o_W^{(**)}$ is coherent.
\end{prop}

Proof. Since $\pi:W^{\rm n}\to W$
is finite, $\pi_*(\o^{\rm n}_W)$ is coherent and torsion free.
Let $w\in W$ be point of codimension $\geq 2$. Then
$\pi^{-1}(w)$ has codimension $\geq 2$ by assumption (2), so
$W^{\rm n}$ has depth $\geq 2$ at $\pi^{-1}(w)$ since  $W^{\rm n} $ is normal.
Thus $\pi_*(\o^{\rm n}_W)$  has depth $\geq 2$ at $w$ and hence
 $\pi_*(\o^{\rm n}_W)$ is $S_2$. Therefore
$$
\o_W^{(**)}\subset \bigl(\pi_*(\o^{\rm n}_W)\bigr)^{(**)}=\pi_*(\o^{\rm n}_W)
\qtq{is coherent.} \qed
$$

Next we show the implications
(\ref{pushfd.finite.thm}.2) $\Leftrightarrow$ (\ref{pushfd.finite.thm}.3) $\Leftrightarrow$ (\ref{pushfd.finite.thm}.4).

\begin{say}[Punctual hulls]\label{punct.hull.say}
Let $(x,X)$ be a  Noetherian, local   scheme and  $F$ a  coherent sheaf on $X$.   Set $U:=X\setminus \{x\}$.  We have an exact sequence
$$
0\to H^0_x(X, F) \to H^0(X, F) \to H^0(U, F|_U)\to H^1_x(X, F)  
\to H^1(X, F)=0.
$$
Thus $H^0(U, F|_U)$ is a finite  $H^0(X, \o_X)$-module iff $ H^1_x(X, F)$ is.
Since $F^{(**)}_x$ is the sheaf corresponding to 
$H^0(U, F|_U)$, we have proved the following.

\medskip
{\it Claim \ref{punct.hull.say}.2.}  $F^{(**)}_x$  is coherent iff
  $H^1_x(X, F)$ has  finite length.\qed
\medskip

Assume next that $F^{(**)}_x$  is coherent and consider its completion
$\widehat{F^{(**)}_x}$. Its depth at $x$ is $\geq 2$ and
$\widehat{F^{(**)}_x}/\hat F\cong F^{(**)}_x/F$ is Artinian.
Thus $\widehat{F^{(**)}_x}$ is the hull of $\hat F$. 
Conversely, if $F^{(**)}_x$  is not coherent, then it is a limit
of an infinite increasing sequence of coherent sheaves
$F\subsetneq F_1\subsetneq  \cdots \subsetneq F^{(**)}_x$
and taking completions shows that the hull of $\hat F$ is also
not coherent. Thus we get the following.

\medskip
{\it Claim \ref{punct.hull.say}.3.}  $F^{(**)}_x$  is coherent iff
$(\hat F)^{(**)}_x$  is  and then
$(\hat F)^{(**)}_x=\widehat{F^{(**)}_x}$.\qed
\medskip

It remains to understand when the punctual hull is coherent
when $(x, X)$ is local and complete. Let $W\subset X$ be an integral subscheme.
Since $X$ is complete, $X$ and hence $W$ 
are excellent (cf.\ \cite[Tag 07QW]{stacks-project}). Thus, by
Proposition \ref{conds.for.coh.hull}, the global hull $\o_W^{(**)}$ is coherent
and so is the punctual hull $(\o_W)_x^{(**)}\subset \o_W^{(**)}$.  

Thus, if
$\hat F$ has no 
1-dimensional associated primes then $F^{(**)}_x$ is coherent
by (\ref{tf.pf.say}.4).
\end{say}

\begin{say}[Modification and localization]  \label{mod.extend.say}
Let $X$ be a  scheme, $F$ a quasi-coherent sheaf on $X$ and $Z\subset  X$ a closed, nowhere dense subscheme. 
A   {\it   modification}  of $F$  {\it centered} at $Z$ is a 
quasi-coherent sheaf $G$ together with a 
map of sheaves
$q:F\to G$ such that
  none of the associated primes of $G$ is contained in $Z$ and  $q$ is an isomorphism over $X\setminus Z$.
A modification is called coherent if $G$ is coherent.
If $\dim Z=0$  then a  modification of $F$ centered at
 $Z$ is called a 
 {\it punctual modification} of $F$.

Let $j:X\setminus Z\into X$ be the natural injection. 
There is a one-to-one
correspondence between  modifications and quasi-coherent sheaves
$$
F/\tors_ZF\subset G\subset j_*\bigl(F|_{X\setminus Z}\bigr).
$$

  Let $\pi: Y\to X$ be a flat morphism and set
$W:=\pi^{-1}(Z)$.
Since $j_*$ commutes with flat base change,
we obtain that
$$
(\pi^*F)^{(**)}_W\cong \pi^*\bigl(F^{(**)}_Z\bigr).
$$
In particular, the local hull is a sheaf in the Zariski topology
and commutes with arbitrary localizations.

Next let $x\in Z$ be a point. By localizing we obtain
$Z_x\subset X_x$ and $F_x$. Let $p_x:F_x\to G_x$ be a 
coherent modification centered at $Z_x$. First we can extend $p_x:F_x\to G_x$
to a coherent modification $p^0:F^0\to G^0$ defined over some
open neighborhood $x\in X^0\subset X$ and then, using  \cite[Exrc.II.5.15]{hartsh},
to a coherent modification $p:F\to G$. The following special case
is  especially useful.
\medskip

{\it Claim \ref{mod.extend.say}.1.}  
Let $X$ be a   Noetherian scheme, $x\in  X$ a point
and $F$ a coherent sheaf on $X$. 
Then every  coherent modification $p_x:F_x\to G_x$  centered at $x\in X_x$ can be extended to a coherent modification $p:F\to G$  centered at $\bar x\in X$.\qed
\end{say}

The following lemma is a special case of the assertion that the definition of
depth  given in  \cite[Exp.III]{sga2}  is equivalent to the usual definition.

\begin{lem}\label{rel.hull.depth.2.lem}
Let $X$ be a   Noetherian scheme and $Z\subset   X$ a closed, nowhere dense subscheme. Let  $F\to G$ be a  coherent modification centered  at $Z$.
Then $G=F^{(**)}_Z$  iff $\depth_ZG\geq 2$.
\end{lem}

Proof. 
We may assume that  $\tors_ZF=\tors_ZG=0$, so
 we have injections  $F\subset G\subset F^{(**)}_Z$ and $G^{(**)}_Z=F^{(**)}_Z$. 
We need to show that $j^{(**)}_Z:G\to G^{(**)}_Z$ is an isomorphism iff
$\depth_ZG\geq 2$. For this we may assume that $X$ is affine.

Pick any $x\in Z$ and let $s\in \o_X$ be an equation of $Z$ that is not a zero divisor on $G$.  Note that  $\depth_xG\geq 2$ iff 
$x$ is not an associated point of $G/sG$. 

If $\depth_xG=1$ then let  $sG\subsetneq G'\subset G$ be a subsheaf such that
$G'/sG$ is supported on $\bar x$. Then  $G\to s^{-1}G'$ is a non-identity
modification of $G$  centered at $\bar x\subset Z$, thus $G\neq G^{(**)}_Z$.

Conversely, if $G\subset G^{(**)}_Z$ is not an isomorphism then there is a  non-identity
coherent modification  $G\to G'$ centered at $Z$. Let $x$ be a generic point of
$\supp (G'/G)$. Then   $s^mG'\subset G$ for some $m>0$ and
$s^mG\subset s^mG'\subset G$ shows that $\depth_xG=1$. \qed

\begin{say}[Proof of Theorem \ref{pushfd.finite.thm}]
 We may assume that $X$ is affine. The equivalence of (\ref{pushfd.finite.thm}.1) and (\ref{pushfd.finite.thm}.5) is 
proved in (\ref{tf.pf.say}.4)
and (\ref{pushfd.finite.thm}.2) $\Leftrightarrow$ (\ref{pushfd.finite.thm}.3) $\Leftrightarrow$ (\ref{pushfd.finite.thm}.4) is proved in
Paragraph \ref{punct.hull.say}.
Note further that the equivalences of (\ref{pushfd.finite.thm}.5), (\ref{pushfd.finite.thm}.6), (\ref{pushfd.finite.thm}.7) and (\ref{pushfd.finite.thm}.8) with each other  are special cases of the equivalences of (\ref{pushfd.finite.thm}.1), (\ref{pushfd.finite.thm}.2), (\ref{pushfd.finite.thm}.3) and (\ref{pushfd.finite.thm}.4) with each other.
Thus it remains to  show that
(\ref{pushfd.finite.thm}.1) and (\ref{pushfd.finite.thm}.2) are equivalent. 

Assume (\ref{pushfd.finite.thm}.1) and set $T=\bar x$. There is a natural map
$ F^{(**)}_T\to F^{(**)}_Z$ whose kernel is supported on $Z$.
Since $W_i\not\subset Z$ for every $i$, $F^{(**)}_T\to F^{(**)}_Z$ is injective.
Therefore $F^{(**)}_T$ is coherent. Since $F^{(**)}_x$
is a localization of $F^{(**)}_T$, 
this shows that (\ref{pushfd.finite.thm}.1) $\Rightarrow$ (\ref{pushfd.finite.thm}.2).

Next assume that (\ref{pushfd.finite.thm}.2) holds. 
By  Proposition \ref{nonS2.finite.inZ.lem} there are only finitely many points $x_i\in Z$ such that $\depth_{x_i}F=1$. Let $D_Z(F)\subset Z$ denote the union of their closures. 

Let $x$ be a generic point of $D_Z(F)$. By our assumption (\ref{pushfd.finite.thm}.2),
$F_x\to F^{(**)}_x$ is a coherent modification  centered at $x$; in particular $\codim_Xx\geq \depth_xF\geq 2$ by Lemma \ref{rel.hull.depth.2.lem}. 
By (\ref{mod.extend.say}.2) we can extend $F^{(**)}_x$  to a coherent modification
$p:F\to F_1$ centered at $\bar x$. 

By (\ref{tf.pf.say}.3) our assumption (\ref{pushfd.finite.thm}.2) also holds for $F_1$ and
we claim that  $D_Z(F_1)\subsetneq D_Z(F)$. To see the latter,
pick any point $x'\in Z$ such that
$\depth_{x'}F_1\leq 1$. If $x'\notin \bar x$ then $p$ is an isomorphism near 
$x'$, hence  $\depth_{x'}F\leq 1$. From Lemma \ref{rel.hull.depth.2.lem} we obtain that  $\depth_{x}F_1\geq 2$,
hence if $x'\in \bar x$ then it is a non-generic point.
Thus  $D_Z(F_1)\subsetneq D_Z(F)$. 

 We can thus
repeat this process and eventually obtain a coherent modification
$F\to F_m$ such that $D_Z(F_m)=\emptyset$. This means that
$\depth_ZF_m\geq 2$ hence  $F_m=F^{(**)}_Z$  by Lemma \ref{rel.hull.depth.2.lem}. 
\qed
\end{say}

\section{Coherence of  global hulls}\label{sec.2}

\begin{say}[Depth of submodules]\label{DEPTH.OF.SUBMODS}
Let $(R, m)$ be a   Noetherian local ring and $M$ a finite
$R$-module. Then  $\depth_mM\geq 1$ iff $M$ has no nonzero submodule
$0\neq N\subset M$ such that $mN=0$.  Thus $\depth_mM'\geq 1$ for every
submodule $M'\subset M$.

Assume next that $\depth_mM\geq 2$ and $M'\subset M$ is $m$-saturated,
meaning that there is no  submodule
$M'\subsetneq N\subset M$ such that $mN\subset M'$.
Thus  there is an $r\in m$ that is 
not a zero divisor on $M/M'$ and so we get an injection
$M'/rM'\into M/rM$. Therefore $\depth_mM'/rM'\geq 1$ and so
$\depth_mM'\geq 2$.

Let next $P\subset R$ be an associated prime of $M$ and
$\tors_PM\subset M$ the torsion submodule corresponding to $P$. 
Then $\tors_PM$ is $m$-saturated, hence
$\dim \tors_PM\geq \depth_m\tors_PM\geq 2$. In particular, 
 every associated prime of $M$ has
dimension $\geq 2$.

Let now $X$ be a  Noetherian  scheme,  $F$  a coherent, $S_2$ sheaf on $X$
and $x\in \supp F$ a point of codimension $\geq 2$. Then
$\depth_xF\geq 2$. Therefore, if $W$ is an irreducible component  of
$\supp F$ that contains $x$ then $\dim \o_{x, W}\geq 2$.
Equivalently,  $x$ also has codimension $\geq 2$ in $W$. 
This is exactly the codimension 1 purity condition
claimed in (\ref{purity.cond.lem}.1).
\end{say}

\begin{prop} \label{nonS2.finite.inZ.lem}
Let $X$ be a Noetherian  scheme, $F$ a coherent sheaf on $X$  and $Z\subset X$ a closed, nowhere dense
subscheme that does not contain any of the associated points of $F$.
Then there are only  finitely many points $x\in Z$ such that
$\depth_xF\leq 1$.
\end{prop}

Proof. The question is local so we may assume that $X$ is affine.
By our assumptions there is  a Cartier divisor $(g=0)$   containing $Z$
and $\depth_xX\leq 1$ iff   $x$ is an associated point of $F/gF$. Since $X$ is Noetherian,
there are only finitely many such points.\qed

\begin{exmp} \label{bad.exmps.2}
The assumption that $Z$ should not contain any of the associated points of $F$ is necessary. For instance, set
 $X=\spec k[x,y,z,t]/(tx, ty, t^2)$.  The associated primes are $(0)$,
$(x,y,t)$ and $\depth_pX=1$ at  every point  $p\in V(x,y,t)$. 

The assumption that $Z$ should be nowhere dense is also necessary.
The following example is modeled on \cite[A.1]{MR0155856}.
Start with $R_1:=k[x_1, y_1, x_2, y_2, \dots]$. Let $R_2\subset R_1$ be the
subring generated by all monomials of degree $\geq 2$ and $R_3$ the ring obtained by inverting every element not contained in any of the ideals
$m_i:=R_2\cap (x_i, y_i)R_1$. Then $R_3$ is Noetherian, has dimension 2 and
its maximal ideals are $m_iR_3$. These form a  Zariski dense subset of $\spec R_3$ and $R_3$ has depth 1 at all of the maximal ideals.
\end{exmp}

The next results says that a generically $S_2$ sheaf is always
$S_2$ outside a codimension $\geq 2$ subset.

\begin{cor} \label{S2.codim.1.ok.cor}
Let $X$ be a Noetherian  scheme, $F$ a coherent, $S_1$  sheaf on $X$  and $Z\subset \supp F$ a closed, nowhere dense
subscheme such that  $F|_{X\setminus Z}$ is $S_2$. Then there is a closed
subscheme $W\subset Z$ such that
\begin{enumerate}
\item $W$ has codimension $\geq 2$ in $\supp F$ and
\item $F|_{X\setminus W}$ is $S_2$.
\end{enumerate}
\end{cor}

Proof. By Proposition \ref{nonS2.finite.inZ.lem},
 there are only  finitely many points $x\in Z$ such that
$\depth_xF\leq 1$. We can thus take  $W:=\cup \bar x$ where
$x$ runs through all points of codimension $\geq 2$ for which
$\depth_xF\leq 1$. \qed

\medskip

\begin{say}[Proof of Theorem \ref{glob.hull.conds.thm}]
The equivalence of (\ref{glob.hull.conds.thm}.1) and 
(\ref{glob.hull.conds.thm}.2) follows from 
(\ref{tf.pf.say}.4).

Assume next that (\ref{glob.hull.conds.thm}.2) holds. 
As we noted in (\ref{purity.cond.lem}.2), if $\o_W$ has a hull then
there is an open dense subset $W^0\subset W$ that is $S_2$ and  
(\ref{pushfd.finite.thm}.5) $\Rightarrow$ (\ref{pushfd.finite.thm}.8)
shows that for every point  $x\in W$ of codimension $\geq 2$, the completion
$\hat{\o}_{x,W}$ has no 1-dimensional associated primes. These are
(\ref{glob.hull.conds.thm}.3.a--b).

Conversely, if (\ref{glob.hull.conds.thm}.3.a) holds then by
Corollary \ref{S2.codim.1.ok.cor} there is a closed
subscheme $W\subset X$ of  codimension $\geq 2$ such that
 $\o_{X\setminus W}$ is $S_2$. Thus  $\o_X^{(**)}=(\o_X)^{(**)}_W$ by
(\ref{loc.glob.tf.exmp}.2) and $(\o_X)^{(**)}_W$ is coherent by 
(\ref{pushfd.finite.thm}.8) $\Rightarrow$ (\ref{pushfd.finite.thm}.5).\qed
\end{say}


\def\cprime{$'$} \def\cprime{$'$} \def\cprime{$'$} \def\cprime{$'$}
  \def\cprime{$'$} \def\cprime{$'$} \def\dbar{\leavevmode\hbox to
  0pt{\hskip.2ex \accent"16\hss}d} \def\cprime{$'$} \def\cprime{$'$}
  \def\polhk#1{\setbox0=\hbox{#1}{\ooalign{\hidewidth
  \lower1.5ex\hbox{`}\hidewidth\crcr\unhbox0}}} \def\cprime{$'$}
  \def\cprime{$'$} \def\cprime{$'$} \def\cprime{$'$}
  \def\polhk#1{\setbox0=\hbox{#1}{\ooalign{\hidewidth
  \lower1.5ex\hbox{`}\hidewidth\crcr\unhbox0}}} \def\cdprime{$''$}
  \def\cprime{$'$} \def\cprime{$'$} \def\cprime{$'$} \def\cprime{$'$}
\providecommand{\bysame}{\leavevmode\hbox to3em{\hrulefill}\thinspace}
\providecommand{\MR}{\relax\ifhmode\unskip\space\fi MR }
\providecommand{\MRhref}[2]{%
  \href{http://www.ams.org/mathscinet-getitem?mr=#1}{#2}
}
\providecommand{\href}[2]{#2}

\bigskip

\noindent  Princeton University, Princeton NJ 08544-1000

{\begin{verbatim} kollar@math.princeton.edu\end{verbatim}}

\end{document}